\documentclass[11pt]{article}

\usepackage[top=50mm, bottom=50mm, left=50mm, right=50mm]{geometry}

\usepackage{lineno}                 
\usepackage{amssymb}                
\usepackage{amsmath}                
\usepackage{amsthm}                 
\usepackage{mathrsfs}               
\usepackage{epsfig}                 
\usepackage{graphicx}               
\usepackage{graphics}               
\usepackage{float}                  %
\usepackage{multirow}               %
\usepackage{color}                  %
\usepackage{fullpage}               %
\usepackage[normalem]{ulem}         %
\usepackage{makeidx}                %
\usepackage{xspace}                 %
\usepackage{wrapfig}                %

\usepackage{url}                    
\usepackage{booktabs}               
\usepackage{amsfonts}               
\usepackage{nicefrac}               
\usepackage{microtype}              

\usepackage{caption}                %
\usepackage{enumitem}               %
\usepackage[export]{adjustbox}      %
\usepackage{comment}                %
\usepackage{xstring}                %
\usepackage{amsfonts}               %
\usepackage{subcaption}             %
\usepackage{tikz}                   %
\usepackage{pgfplots}               %
\usepackage{pgfplotstable}          %

\usepackage[affil-it]{authblk}                

\usepackage[citestyle=numeric-comp,sorting=none,backend=biber]{biblatex}
\makeindex

\providecommand{\keywords}[1]
{
  \small	
  \textbf{\textit{Keywords: }} #1
}

\pgfplotsset{compat=1.17}

\newtheorem{problem}{Problem}

\newtheorem{remark}{Remark}

\newcommand\mydef{:=}

\def\bsigma{\hbox{\boldmath$\sigma$}}
\def\bepsilon{\hbox{\boldmath$\epsilon$}}
\def\bnabla{\hbox{\boldmath$\nabla$}}

\def\disp{\hbox{$u$}}
\def\bdisp{\hbox{$\mathbf{u}$}}

\def\testbdisp{\hbox{$\delta\mathbf{u}$}}

\def\strain{\hbox{$\bepsilon$}}

\def\stress{\hbox{$\bsigma$}}
\def\stresspos{\hbox{$\bsigma^+$}}
\def\stressneg{\hbox{$\bsigma^-$}}
\newcommand{\lame}[1]{%
    \IfEqCase{#1}{%
        {1}{\lambda}%
        {2}{\mu}%
    }[\PackageError{lame}{Undefined option to lame: #1}{}]%
}

\def\pf{\hbox{$\varphi$}}
\def\testpf{\hbox{$\delta\varphi$}}
\def\l{\hbox{$l$}}
\def\gc{\hbox{$G_c$}}

\def\slk{\hbox{$\theta$}}

\def\dom{\hbox{$\Omega$}}
\def\surf{\hbox{$\Gamma$}}

\newcommand{\surfarg}[2]{\hbox{$\surf_{#1}^{\scriptsize{#2}}$}}
\newcommand{\trac}[1]{\hbox{$\mathbf{t}_{p}^{\scriptsize{#1}}$}}

\usetikzlibrary{fadings}
\tikzfading[name=fade out, 
    inner color=transparent!0,
    outer color=transparent!100]
\usetikzlibrary{shapes}
\usetikzlibrary{decorations.pathreplacing}
\usetikzlibrary{arrows}
\usepgfplotslibrary{fillbetween}

\DeclareCaptionFont{mysize}{\fontsize{9}{9.6}\selectfont}
\begin{document}

\title{Phase-field fracture irreversibility using the slack variable approach}

\author[1]{Ritukesh Bharali}
\author[1]{Fredrik Larsson}
\author[2]{Ralf J\"anicke}
\affil[1]{Department of Industrial and Material Science, Chalmers University of Technology}
\affil[2]{Institute of Applied Mechanics, Technische Universit\"at Braunschweig}

\date{}


\maketitle

\section*{Highlights}
\begin{itemize}
    \item Effective equality constraint for fracture irreversibility inequality constraint using slack variable.
    \item The Lagrange Multiplier Method and the Penalty method are adopted to augment the energy functional.
    \item Numerical experiments carried out on both brittle and quasi-brittle fracture problems. 
\end{itemize}

\bigskip

\begin{abstract}
\noindent In this manuscript, the phase-field fracture irreversibility constraint is transformed into an equality-based constraint using the slack variable approach. The equality-based fracture irreversibility constraint is then introduced in the phase-field fracture energy functional using the Lagrange Multiplier Method and the Penalty method. Both methods are variationally consistent with the conventional variational inequality phase-field fracture problem, unlike the history-variable approach. Thereafter, numerical experiments are carried out on benchmark problems in brittle and quasi-brittle fracture to demonstrate the efficacy of the proposed method.
\end{abstract}

\keywords{phase-field fracture, brittle, quasi-brittle, Lagrange multiplier method, Penalty method, COMSOL}

\newpage

\section{Introduction}

The phase-field model for fracture emerged from the \textit{seminal work} of \cite{Francfort1998}, wherein the Griffith fracture criterion was cast into a variational setting. Later, a numerical implementation of the same was proposed in \cite{Bourdin2000}, using the Ambrosio--Tortorelli regularisation of the Mumford-Shah potential \cite{Mumford1989}. In this implementation, the fracture is represented by an auxiliary variable, that interpolates between the intact and broken material states. Such a formulation allows the automatic tracking of fractures on a fixed mesh, thereby eliminating the need for the tedious tracking and remeshing processes observed with discrete methods. Furthermore, the phase-field model for fracture is also able to handle topologically complex (branching, kinking and merging) fractures, and is able to demonstrate fracture initiation without introducing any singularity \cite{GERASIMOV2019990}. Owing to these advantages, the phase-field model for fracture has gained popularity in the computational mechanics community in the past decade.

A thermodynamically consistent formulation of the phase-field fracture model was proposed by \cite{Miehe2010a}, adopting an energetic cracking driving force definition. Since then, over the past decade, researchers have extended the work to thermo-mechanical problems at large strain \cite{Miehe2015a}, ductile failure \cite{Miehe2015b,Ambati2015a}, fracture in thin films \cite{Mesgarnejad2013}, anisotropic fracture \cite{Nguyen2017}, fracture in fully/partially saturated porous media \cite{Miehe2016b,Cajuhi2018,Mikelic2019}, hydrogen assisted cracking \cite{Martinez-Paneda2018}, dissolution-driven stress corrosion cracking \cite{CUI2021104254}, fracture and fatigue in shape-memory alloys \cite{SIMOES2021113504}, and brittle failure of Reisner-Mindlin plates \cite{KIKIS2021113490} to cite a few. Finally, the phase-field model for fracture has also been used for the multi-scale finite element method \cite{Patil2018a,Patil2018b,Patil2019}, asymptotic homogenisation \cite{Fantoni2019}, and variationally consistent homogenisation \cite{BHARALI2021104247}. For a detailed overview on the phase-field fracture model, the reader is referred to the comprehensive review works \cite{Ambati2014,WU20201,de2020numerical}.

Despite its growing popularity, the phase-field model for fracture poses several challenges when it comes to robust and computationally efficient solution techniques \cite{de2020numerical}. These include,

\begin{enumerate}[label=(\Alph*)]
    \item the poor performance of the monolithic Newton-Raphson (NR) method due to the non-convex energy functional,
    \item variational inequality arising from fracture irreversibility, and 
    \item extremely fine meshes are required in the fracture zone.
\end{enumerate}

In order to alleviate the problem concerning the poor-performance of the NR method (A), \cite{Heister2015} utilised a \textit{linear extrapolation in time}\footnote{The linear extrapolation of the phase-field variable in time is a questionable assumption although it results in a robust solution method.} for the phase-field variable in the momentum balance equation. In \cite{Gerasimov2016}, a novel line search technique was developed, while in \cite{Wick2017a} and \cite{Wick2017b} Augmented Lagrangian and modified NR approaches were proposed. More recently, \cite{KOPANICAKOVA2020112720} adopted a recursive multilevel trust region method, that resulted in improved convergence of the NR method. The development of robust, monolithic solution technique is still an active area in phase-field fracture research, and this manuscript is contribution towards this aspect. 

The next issue pertains to the variational inequality problem arising out of the fracture irreversibility condition (B). This aspect has been treated in different ways in the phase-field fracture literature. They include crack set irreversibility \cite{Bourdin2000,Bourdin2007,Burke2010}, penalisation \cite{GERASIMOV2019990,Wick2017a,Wick2017b}, and the implicit History variable based method proposed by \cite{Miehe2010}. The lattermost method remains popular despite its non-variational nature, over-estimation of the fracture surface energy, and the necessity for computationally expensive alternate minimisation (staggered) solution procedure \cite{GERASIMOV2019990}. In a novel approach, this manuscript adopts the slack variable method \cite{valentine1933problem,Bertsekas1997}, that constructs an equivalent variational equality problem, while maintaining the variational structure of the original problem.

Another issue pertaining to the phase-field fracture model is the requirement of extremely fine meshes in the smeared fracture zone (C). In this context, \cite{Bourdin2000} and \cite{Mesgarnejad2013} advocates the use of uniformly refined meshes together with parallel computing while \cite{burke2010adaptive,artina2015anisotropic,wick2016goal} opted for error-controlled adaptive mesh refinement. Some other approaches include pre-refined meshes when the fracture path is known \cite{Ambati2014,Gerasimov2016,GERASIMOV2019990}, adaptive mesh refinement in \cite{Heister2015} based on a phase-field threshold or Kelly error estimates \cite{Kelly1983}, and multi-level hp-FEM in \cite{Nagaraja2019}. In this manuscript, uniformly refined meshes and pre-refined meshes based on a priori knowledge of the fracture path, are used.

The novelty of this manuscript lies in the alternative treatment of the fracture irreversibility inequality constraint, using the slack variable approach. The inequality constraint is replaced by an equivalent equality-based fracture irreversibility constraint. The constraint is then augmented to the phase-field fracture energy functional using the Lagrange Multiplier method and the Penalty method. These methods are variationally consistent, unlike the history-variable approach proposed in \cite{Miehe2010}.

The manuscript is organised as follows: Section \ref{sec2} introduces the phase-field model for fracture, its underlying energy functional and the Karush-Kuhn-Tucker (KKT) optimality conditions. Thereafter, in Section \ref{sec3}, the Method of Multipliers is adopted to augment the energy functional a slack variable-based fracture irreversibility criterion. The equivalence of the new formulation with the original problem in Section \ref{sec2} is established in terms of KKT conditions. Furthermore, the Euler-Lagrange equations are presented. In Section \ref{sec4}, the slack variable-based fracture irreversibility constraint is introduced in the energy functional using the Penalisation method, followed by the derivation of the pertinent Euler-Lagrange equations. The numerical experiments on benchmark problems are presented in Section \ref{sec5}, and Section \ref{sec6} lays down the concluding remarks of this manuscript.

\section*{Notation}

\begin{em}
The following notations are strictly adhered to in this manuscript:
\begin{itemize}
\item Zero-order tensors are represented using small italicized letters, e.g., $a$. Bold symbols are used for first and second-order tensors, for instance, stress $\stress$ and strain $\strain$.
\item A function $f$ with its arguments $x$, $y$ is written in the form $f(x,y)$, whereas a variable $g$ with operational dependencies $p,q$ is written as $g[p,q]$.
\item The Macaulay operator on a variable $x$ is defined as $\langle x \rangle_\pm = \frac{1}{2} (x \pm |x|)$.
\end{itemize}
\end{em}

\section{Phase-field fracture model}\label{sec2}

\subsection{Energy functional}

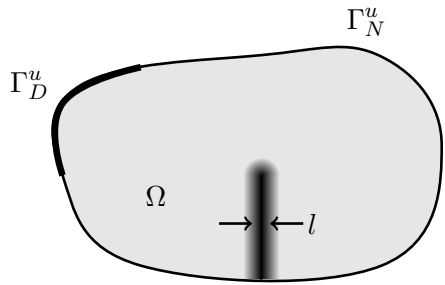
\begin{figure}[ht!]
    \centering
    \begin{tikzpicture}[scale=0.8]
    \coordinate (K) at (0,0);
    \draw [fill=black!10,line width=1pt] (K) plot [smooth cycle,tension=0.7] 
    coordinates {(3,1) (7,1) (8,3) (7,4.475) (5,4.5) (2,4) (1.7,2.5)};
    \node[ ] at (3.25,2.15) {$\dom$};
    \fill[black, path fading=fade out, draw=none] (5,2.5) circle (0.3);
    \draw[line width=0.1pt,black] (5,0.75) to (5,2.5);
    \draw[line width=0.1pt,black] (5,0.75) to (5,2.5);
    \shade [top color=black,bottom color=black!10,shading angle=90] (5,0.78) rectangle (5.3,2.5);
    \shade [top color=black!10,bottom color=black,shading angle=90] (4.7,0.78) rectangle (5,2.5);
    \draw[->,line width=1pt,black] (4.3,1.7) to (4.85,1.7);
    \draw[<-,line width=1pt,black] (5.15,1.7) to (5.7,1.7);
    \node[ ] at (5.85,1.7) {$\l$};
    \draw (K) [line width=2.5pt,black] plot [smooth, tension=0.8] coordinates {(3,4.3) (1.7,3.7) (1.7,2.5)};
    \node[ ] at (1.15,4) {$\surfarg{D}{u}$};
    \node[ ] at (6.75,5.05) {$\surfarg{N}{u}$};
    \end{tikzpicture}
  \caption{A solid $\dom \in \mathbb{R}^2$ embedded with a diffused (smeared) crack. Dirichlet and Neumann boundaries are indicated with $\surfarg{D}{u}$ and $\surfarg{N}{u}$ respectively. Figure adopted from \cite{BHARALI2021104247}.}
    \label{fig:sec2:continuumpotato}
\end{figure}

Let $\dom \in \mathbb{R}^{\text{dim}}$ ($\text{dim} = 2,3$) be the domain occupied by the fracturing solid, shown in Figure \ref{fig:sec2:continuumpotato}. Its boundary $\surf$ is decomposed into a Dirichlet boundary $\surfarg{D}{u}$ and a Neumann boundary $\surfarg{N}{u}$, such that $\surf = \surfarg{D}{u} \cup \surfarg{N}{u}$ and $\surfarg{D}{u} \cap \surfarg{N}{u} = \emptyset$. Furthermore, the fracture is represented by an auxiliary variable (phase-field) $\pf \in [0,1]$ within a diffusive (smeared) zone of width $\l>0$.

\noindent The energy functional of the phase-field fracture model \cite{Gerasimov2016,GERASIMOV2019990} is stated as,

\begin{equation}{\label{eqn:sec2:EFunc}}
    \displaystyle E(\bdisp,\pf) = \int_{\dom}^{} g(\pf) \Psi^+(\strain[\bdisp]) \: \text{d}\dom + \int_{\dom}^{} \Psi^-(\strain[\bdisp]) \: \text{d}\dom - \int_{\surfarg{N}{u}}^{}  \trac{u} \cdot \bdisp \: \text{d}\surf + \int_{\dom}^{} \dfrac{\gc}{c_w \l } \left( w(\pf) + \l |\bnabla \pf|^2 \right) \: \text{d}\dom,
\end{equation}

\noindent with $\bdisp$ and $\pf$ representing the displacement and the phase-field respectively. The symbols adopted for the additional functions and parameters in (\ref{eqn:sec2:EFunc}) are presented in Table \ref{sec2:table:eFuncParams}.

\begin{table}[ht!]
    \centering
    \begin{tabular}{ll} \hline
    Symbol  & Function/Parameter  \\ \hline
    $g(\pf)$  & Degradation function \\ 
    $\Psi^+$  & Fracture driving strain energy density \\ 
    $\Psi^-$  & Residual strain energy density \\
    $\strain[\bdisp]$  & Symmetric gradient of $\bdisp$ \\
    $\trac{u}$  & Prescribed traction \\ 
    $\gc$  & Griffith fracture energy \\
    $c_w$  & Normalisation constant \\ 
    $\l$  & Fracture length-scale \\
    $w(\pf)$  & Local fracture energy function \\ \hline
    \end{tabular}
    \caption{List of functions and parameters corresponding to the phase-field fracture energy functional (\ref{eqn:sec2:EFunc}).}
    \label{sec2:table:eFuncParams}
\end{table}

In this manuscript, $\Psi^+$ and $\Psi^-$ represent the tensile and the compressive strain energy density functions respectively. They are functions of strain tensor $\strain[\bdisp]$, given by

\begin{equation}{\label{eqn:sec2:def_strain}}
    \strain[\bdisp] = (\bdisp \otimes \bnabla)^{\text{sym}}.
\end{equation}

\noindent The explicit expression for $\Psi^{\pm}$ is obtained following \cite{Miehe2010a} as,

\begin{equation}{\label{eqn:sec2:defPsipm}}
    \Psi^\pm(\strain[\bdisp]) = \frac{1}{2} \lame{1} \langle (\mathbf{I} \colon \strain[\bdisp]) \rangle^2_\pm + \lame{2} \strain^\pm[\bdisp] \colon \strain^\pm[\bdisp],
\end{equation}

\noindent where $\lame{1}$ and $\lame{2}$ are the Lam\'e constants, $\mathbf{I}$ is a second-order identity tensor, $\bdisp$ is the displacement. Moreover, the tensile/compressive strain $\strain^\pm[\bdisp]$ is defined as 

\begin{equation}\label{eqn:sec2:strain_split}
	\strain^\pm[\bdisp] = \sum_{i = 1}^{\text{dim}} \langle \epsilon_i \rangle_\pm \mathbf{p}_i \otimes \mathbf{p}_i
\end{equation} 

\noindent where $\epsilon_i$ represents the $i^{th}$ eigenvalue of the strain, and $\mathbf{p}_i$ is its corresponding normalised eigenvector. The tensile/compressive strain energy definition in (\ref{eqn:sec2:defPsipm}) yield the corresponding stresses upon taking the derivative w.r.t the strain. They are explicitly stated in a compact form as,

\begin{equation}\label{sec2:stress_def}
	\stress^\pm = \frac{\partial \Psi^\pm}{\partial \strain} = \lame{1} \langle (\mathbf{I} \colon \strain[\bdisp]) \rangle_\pm \mathbf{I} + 2 \lame{2} \strain^\pm[\bdisp].
\end{equation}

The last integral in the energy functional (\ref{eqn:sec2:EFunc}) represents the fracture energy in a phase-field fracture model. The phase-field fracture model allows flexibility in the choice of the degradation function $g(\pf)$ and the locally dissipated fracture energy function $w(\pf)$. Table \ref{sec2:table:pf_models} presents commonly adopted functions for brittle and quasi-brittle fracture.

\begin{table}[ht!]
    \centering
    \begin{tabular}{llll} \hline
     Model        & $g(\pf)$  & $w(\pf)$  & $c_w$ \\ \hline
    Brittle-AT1\cite{pham2011}   & $(1-\pf)^2$  & $\pf$  & $8/3$ \\ 
    Brittle-AT2\cite{Bourdin2007}   & $(1-\pf)^2$  & $\pf^2$  & $2$ \\ 
    Quasi-Brittle\cite{wu2017}  & $\dfrac{(1-\pf)^p}{(1-\pf)^p + a_1\pf + a_1 a_2 \pf^2 + a_1 a_2 a_3 \pf^3}$  & $2\pf-\pf^2$  & $\pi$ \\ \hline
    \end{tabular}
    \caption{Degradation function, local fracture energy function and its normalisation constant for different phase-field fracture models. The abbreviation `AT' stands for Ambrosio-Tortorelli.}
    \label{sec2:table:pf_models}
\end{table}

\noindent For the Quasi-Brittle model in Table \ref{sec2:table:pf_models}, the parameters $p$, $a_1$, $a_2$ and $a_3$ are chosen such that they mimic the different traction-separation laws in cohesive zone modelling. Following \cite{wu2017}, the constant $a_1$ is given by

\begin{equation}
    a_1 = \dfrac{4 E_0 \gc}{\pi \l f_t^2}
\end{equation}

\noindent where, the newly introduced parameters $E_0$ and $f_t$ represent the Youngs' modulus and the tensile strength of the material. Next, Table \ref{sec2:table:quasiBrittleParams} presents the different traction-separation laws and the corresponding values of $p$, $a_2$ and $a_3$ \cite{wu2017}.

\begin{table}[ht!]
    \centering
    \begin{tabular}{llll} \hline
     Model        & $p$  & $a_2$  & $a_3$ \\ \hline
    Linear Softening   & $2$  & $-0.5$  & $0$ \\ 
    Exponential Softening   & $2.5$  & $2^{5/3}-3$  & $0$ \\ 
    Cornelissen et. al.\cite{cornelissen1986experimental} Softening  & $2$  & $1.3868$  & $0.6567$ \\ \hline
    \end{tabular}
    \caption{Quasi-brittle phase-field fracture models and their parameters, based on \cite{wu2017}.}
    \label{sec2:table:quasiBrittleParams}
\end{table}

The prediction of fracture initiation and propagation in the solid occupying the domain $\dom$ requires solving the constrained minimisation problem:

\begin{problem}\label{Problem1}
Find $\bdisp(t)$ and $\pf(t)$, where ${}^{n}\pf = \pf(t_n)$ is known for $t_n<t$ such that,
\begin{equation}\label{eqn:sec2:min_E}
    \begin{array}{lll}
    min_{\tiny{\bdisp(t),\pf(t)}} E  & \text{and} & h(\pf) \geq 0,
\end{array}
\end{equation}

\noindent where $h(\pf) = \pf - {}^{n}\pf$, and the left superscript $n$ denotes the previous time-step. Moreover, the problem is augmented with (pseudo) time-dependent Dirichlet boundary conditions $\bdisp^p$ on $\surfarg{D}{u}$ and $\pf^p$ on $\surfarg{D}{\pf}$, and/or Neumann boundary conditions $\trac{u}$ on $\surfarg{N}{u}$ and $\trac{\pf}$ on $\surfarg{N}{\pf}$. Furthermore, the boundary $\surf$ is split as $\surf = \surfarg{D}{u} \cup \surfarg{N}{u}$, $\surfarg{D}{u} \cap \surfarg{N}{u} = \emptyset$ and $\surf = \surfarg{D}{\pf} \cup \surfarg{N}{\pf}$, $\surfarg{D}{\pf} \cap \surfarg{N}{\pf} = \emptyset$ respectively. {\color{black}\hfill $\blacksquare$}
\end{problem}

\subsection{Karush-Kuhn-Tucker conditions}

The KKT triple ($\bdisp,\pf,\Lambda$) in appropriate spaces for the Problem \ref{Problem1} requires

\begin{subequations}
\begin{align}
\big[E(\bdisp,\pf) - \int_{\dom} \Lambda h(\pf) \, \text{d}\dom \big]'_{\bdisp}  & = 0, \label{eqn:sec2:stationarity1}  \\
\big[E(\bdisp,\pf) - \int_{\dom} \Lambda h(\pf) \, \text{d}\dom \big]'_{\pf}  & = 0,, \label{eqn:sec2:stationarity2}  \\
h(\pf) & \geq 0, \label{eqn:sec2:primal_feasibility} \\
\Lambda & \geq 0, \label{eqn:sec2:dual_feasibility} \\
\Lambda h(\pf) & = 0 \label{eqn:sec2:complementary_slackness}
\end{align}
\end{subequations}

\noindent to hold, where $\Lambda$ is a multiplier. Equations (\ref{eqn:sec2:stationarity1}) and (\ref{eqn:sec2:stationarity2}) represent the Euler-Lagrange equations pertaining to the stationary condition, while equations (\ref{eqn:sec2:primal_feasibility}-\ref{eqn:sec2:complementary_slackness}) are the primal and dual feasibility conditions, and the complementary slackness respectively. Formally, the latter set of equations enforces fracture irreversbility. To the same end, alternative, but equivalent formulations have been adopted in the phase-field fracture literature. For instance, \cite{Heister2015} adopted the primal-dual active set strategy proposed by \cite{Hintermuller2002}, \cite{Wick2017a} utilised an Augmented Lagrangian penalisation approach with the Moreau-Yoshida indicator function, and \cite{GERASIMOV2019990} opted for a simple penalisation approach. The authors would like to emphasise that the rather popular history-variable approach in \cite{Miehe2010a} is not variationally consistent.

\section{Lagrange Multiplier Method (LMM)}\label{sec3}

\subsection{Fracture irreversibility and modified energy functional}

In order to enforce fracture irreversibility $h(\pf) = \pf - {}^{n}\pf \geq 0$, a slack variable is defined as 

\begin{equation}\label{eqn:sec3:slk_def}
\slk^2 = h(\pf) = \pf - {}^{n}\pf.  
\end{equation}

\noindent It is observed that $\slk^2$ admits value greater than or equal to zero, thereby fulfilling the fracture irreversibility criterion. Next, the constrained minimisation Problem \ref{Problem1} is reformulated as:

\begin{problem}\label{Problem2}
Find $\bdisp(t)$, $\pf(t)$ and $\slk(t)$ such that,
\begin{equation}\label{eqn:sec3:min_E}
    \begin{array}{lll}
    \text{min}_{\scriptsize{\bdisp(t),\pf(t)},\slk(t)} \:\: E \\
    \text{subjected to} \:\: h(\pf) - \slk^2 = 0. 
\end{array}
\end{equation}

\noindent with suitable (pseudo) time-dependent Dirichlet and/or Neumann boundary conditions, as mentioned in Problem \ref{Problem1}. {\color{black}\hfill $\blacksquare$}
\end{problem}

\noindent Augmenting the equality constraint (\ref{eqn:sec3:min_E}) in (\ref{eqn:sec2:EFunc}) via the Lagrange multiplier $\Lambda$ results in the modified energy functional,

\begin{equation}{\label{eqn:sec3:Efunc_slk}}
\Tilde{E}(\bdisp,\pf,\theta,\Lambda) = E(\bdisp,\pf) - \Lambda \big( h(\pf)-\slk^2 \big).    
\end{equation}

\subsection{Stationary conditions}

The stationary conditions for the Problem \ref{Problem2} are given by

\begin{subequations}
\begin{align}
\Tilde{E}'_{\bdisp}(\bdisp,\pf,\slk,\Lambda; \testbdisp)  & = 0, \label{eqn:sec3:stationarity1} \\
\Tilde{E}'_{\pf}(\bdisp,\pf,\slk,\Lambda; \testpf)  & = 0, \label{eqn:sec3:stationarity2} \\
\Tilde{E}'_{\slk}(\bdisp,\pf,\slk,\Lambda; \delta\slk)  & = 0, \label{eqn:sec3:stationarity3} \\
\Tilde{E}'_{\Lambda}(\bdisp,\pf,\slk,\Lambda; \delta\Lambda)  & = 0. \label{eqn:sec3:stationarity4}
\end{align}
\end{subequations}

\noindent For equivalence of these conditions with those presented for Problem \ref{Problem1}, the reader is referred to Proposition 3.1 in \cite{Tapia1974}. 


\subsection{Euler-Lagrange equations}

The Euler-Lagrange equations for Problem \ref{Problem2} are (\ref{eqn:sec3:stationarity1})-(\ref{eqn:sec3:stationarity4}), obtained upon taking the first variation of the energy functional (\ref{eqn:sec3:Efunc_slk}) w.r.t its solution variables $\bdisp$, $\pf$, $\theta$ and $\Lambda$. This results in,

\begin{problem}\label{Problem3}
Find ($\bdisp$, $\pf$, $\theta$, $\Lambda$) $\in \mathbb{U} \times \mathbb{V} \times \mathbb{W} \times \mathbb{A}$ with

\begin{subequations}
\begin{align}
R^{\mathbf{u}} & = \int_{\dom}^{} \bigg( g(\pf) \stresspos[\bdisp] + \stressneg[\bdisp] \bigg) \colon \strain[\testbdisp] \: \normalfont\text{d}\dom = 0 \:\:\: & \forall \: \testbdisp \in \mathbb{U}^0 \\
R^{\varphi} & = \int_{\dom}^{} \bigg( \dfrac{\gc}{c_w\l} w'(\pf) + g'(\pf) \Psi^+(\strain[\bdisp]) - \Lambda \bigg) \testpf \: \normalfont\text{d}\dom + \int_{\dom}^{} \dfrac{2\gc\l}{c_w} \bnabla\pf \cdot \bnabla\testpf \: \normalfont\text{d}\dom = 0 \:\:\: & \forall \: \testpf \in \mathbb{V}^0 \\
R^{\theta} & = \int_{\dom}^{} 2 \Lambda \theta \delta\theta \: \normalfont\text{d}\dom = 0 \:\:\: & \forall \: \delta\theta \in \mathbb{W} \\
R^{\Lambda} & = - \int_{\dom}^{} \big( h(\pf)-\slk^2 \big) \delta\Lambda \: \normalfont\text{d}\dom = 0 \:\:\: & \forall \: \delta\Lambda \in \mathbb{A}.
\end{align}
\end{subequations}

with trial function spaces
\begin{subequations}
\begin{align}
&\mathbb{U} \mydef \{ \bdisp \in [H^1(\dom)]^{\normalfont \text{dim}} | \bdisp = \bdisp^{\text{{\normalfont p}}} \:\: \text{\normalfont on} \:\: \surfarg{D}{u} \}, \\ 
&\mathbb{V} \mydef \{ \pf \in [H^1(\dom)]^1 | \pf = \pf^{\text{{\normalfont p}}} \:\: \text{\normalfont on} \:\: \surfarg{D}{\pf} \}, \\
&\mathbb{W} \mydef \{ \theta \in [L_2(\dom)] \}, \\
&\mathbb{A} \mydef \{ \Lambda \in [L_2(\dom)] \},
\end{align}
\end{subequations}

and  test function spaces

\begin{subequations}
\begin{align}
&\mathbb{U}^0 \mydef \{ \bdisp \in [H^1(\dom)]^{\normalfont \text{dim}} | \bdisp = 0 \:\: \text{\normalfont on} \:\: \surfarg{D}{u} \}, \\ 
&\mathbb{V}^0 \mydef \{ \pf \in [H^1(\dom)]^1 | \pf = 0 \:\: \text{\normalfont on} \:\: \surfarg{D}{\pf} \},
\end{align}
\end{subequations}

and pertinent (pseudo) time-dependent Neumann condition

\begin{equation}
\mathbf{t} \mydef \left( g(\pf) \stresspos + \stressneg \right) \cdot\mathbf{n} = \trac{u} \:\: \text{\normalfont on} \:\: \surfarg{N}{u}. \label{eqn:sec3:trac_def} 
\end{equation}
{\color{black}\hfill $\blacksquare$}
\end{problem}

\begin{remark}
In the event $\Lambda$ and $\slk$ are zero, their stiffness contribution is zero and the stiffness matrix is singular. In order to avoid this situation a fictitious stiffness is introduced to guarantee solvability. 
\end{remark}

\section{Penalty (Pen.) method}\label{sec4}

\subsection{The energy functional}

In an alternative approach, the squared slack variable equation (\ref{eqn:sec3:slk_def}) is introduced in the energy functional (\ref{eqn:sec2:EFunc}) as a quadratic term via penalisation. This results in the modified energy functional,

\begin{equation}{\label{eqn:sec4:Efunc_slk}}
\hat{E}(\bdisp,\pf,\theta) = E(\bdisp,\pf) + \dfrac{\eta}{2} \int_{\dom} \bigg(h(\pf)-\slk^2\bigg)^2 \, \text{d}\dom,   
\end{equation}

\noindent where, $\eta$ is a penalty parameter.

\subsection{Euler-Lagrange equations}

The Euler-Lagrange equations are obtained upon taking the first variation of the energy functional (\ref{eqn:sec4:Efunc_slk}) w.r.t its solution variables $\bdisp$, $\pf$ and $\theta$. This results in:

\begin{problem}\label{Problem4}
Find ($\bdisp$, $\pf$, $\theta$) $\in \mathbb{U} \times \mathbb{V} \times \mathbb{W}$ with

\begin{subequations}
\begin{align}
R^{\mathbf{u}} & = \int_{\dom}^{} \bigg( g(\pf) \stresspos[\bdisp] + \stressneg[\bdisp] \bigg) \colon \strain[\testbdisp] \: \normalfont\text{d}\dom = 0 \:\:\: & \forall \: \testbdisp \in \mathbb{U}^0 \\
R^{\varphi} & = \int_{\dom}^{} \bigg( \dfrac{\gc}{c_w\l} w'(\pf) + g'(\pf) \Psi^+(\strain[\bdisp]) + \eta (h(\pf)-\slk^2) \bigg) \testpf \: \normalfont\text{d}\dom + \int_{\dom}^{} \dfrac{2\gc\l}{c_w} \bnabla\pf \cdot \bnabla\testpf \: \normalfont\text{d}\dom = 0 \:\:\: & \forall \: \testpf \in \mathbb{V}^0 \\
R^{\theta} & = \int_{\dom}^{} -2 \eta \theta (h(\pf)-\slk^2) \delta\theta \: \normalfont\text{d}\dom = 0 \:\:\: & \forall \: \delta\theta \in \mathbb{W}
\end{align}
\end{subequations}

with trial function spaces
\begin{subequations}
\begin{align}
&\mathbb{U} \mydef \{ \bdisp \in [H^1(\dom)]^{\normalfont \text{dim}} | \bdisp = \bdisp^{\text{{\normalfont p}}} \:\: \text{\normalfont on} \:\: \surfarg{D}{u} \}, \\ 
&\mathbb{V} \mydef \{ \pf \in [H^1(\dom)]^1 | \pf = \pf^{\text{{\normalfont p}}} \:\: \text{\normalfont on} \:\: \surfarg{D}{\pf} \}, \\
&\mathbb{W} \mydef \{ \theta \in [L_2(\dom)]| \},
\end{align}
\end{subequations}

and  test function spaces

\begin{subequations}
\begin{align}
&\mathbb{U}^0 \mydef \{ \bdisp \in [H^1(\dom)]^{\normalfont \text{dim}} | \bdisp = 0 \:\: \text{\normalfont on} \:\: \surfarg{D}{u} \}, \\ 
&\mathbb{V}^0 \mydef \{ \pf \in [H^1(\dom)]^1 | \pf = 0 \:\: \text{\normalfont on} \:\: \surfarg{D}{\pf} \}, 
\end{align}
\end{subequations}

and pertinent (pseudo) time-dependent Neumann condition

\begin{equation}
\mathbf{t} \mydef \left( g(\pf) \stresspos + \stressneg \right) \cdot\mathbf{n} = \trac{u} \:\: \text{\normalfont on} \:\: \surfarg{N}{u}. \label{eqn:sec4:trac_def} 
\end{equation}
{\color{black}\hfill $\blacksquare$}
\end{problem}

\begin{remark}
In the event $h(\pf)-\slk^2$ and $\slk$ are zero, their stiffness contribution is zero and the stiffness matrix is singular.  In order to avoid this situation a fictitious stiffness is introduced to guarantee solvability. 
\end{remark}

\section{Numerical Study}\label{sec5}

The numerical experiments on benchmark problems are presented in this section. They include both brittle and quasi-brittle fracture. In the brittle fracture domain, the Single Edge Notched Specimen under Tension (SENT), and under Shear (SENS), and the Notched concrete specimen with a hole \cite{Ambati2014} is presented. For quasi-brittle fracture, the benchmark problems include the concrete three-point bending problem \cite{Rots1988}, and the Winkler L-shaped panel experiment \cite{winkler2001traglastuntersuchungen}. For each problem, the geometry as well as the material properties and loading conditions are presented in the respective sub-section. The load-displacement curves and the phase-field topology at the final step of analysis is also presented.


For the penalty method, the penalty parameter $\eta$ is set to $10^6$ [N/mm$^2$] throughout the study. This choice is motivated by similarity of the load-displacement curves obtained with those from the Lagrange Multiplier Method. Furthermore, for both the aforementioned methods, a fully coupled monolithic solution technique is adopted using the Newton-Raphson solver in COMSOL Multiphysics. The iterations are terminated when the error measure defined as the weighted Euclidean norm of the solution update,

\begin{equation}
    err: = \sqrt{\dfrac{1}{M}} \sqrt{ \sum_{j=1}^{M} \dfrac{1}{N_j} \sum_{i=1}^{N} \bigg( \dfrac{|E_{i,j}|}{W_{i,j}} \bigg)^2 },
\end{equation}

\noindent is less than $1e-4$. In the above equation, $M$ indicates the number of fields, $j$ is the number of degree of freedom for each field $j=1,2,..M$, $E$ represents the absolute update (say, $\bdisp^{(m)}$-${\bdisp}^{(m-1)}$ for displacement, $m$ being the iteration count), and $W_{i,j} = max(|\mathbf{U}_{i,j}|,S_j)$. The entire solution vector is represented with $\mathbf{U}$ and $S_j$ refers to scaling of solution variables\footnote{Scaling of solution variables prevents possible ill-conditioning of the stiffness matrix.}. 

\subsection{Brittle: Single Edge Notched specimen under Tension (SENT)}\label{sec5:SENtension}

A unit square (in mm) embedded with a horizontal notch, midway along the height is considered, as shown in Figure \ref{sec4:fig:tension}. The length of the notch (shown in red) is equal to half of the edge length of the plate. The notch is modelled explicitly in the finite element mesh. A quasi-static loading is applied at the top boundary in the form of prescribed displacement increment $\Delta\disp = 1e-5$[mm] for the first $450$ steps, following which it is changed to $1e-6$[mm]. Furthermore, the bottom boundary remains fixed. The additional model parameters are presented in Table \ref{sec4:table:miehe_tension_shear}.

\begin{figure}[ht]
\begin{minipage}[b]{0.33\linewidth}
\centering
\begin{tikzpicture}[scale=0.7]
    \coordinate (K) at (0,0);
    \draw[line width=0.75pt,black] (-2.5,-2.5) to (2.5,-2.5);
    \draw[line width=0.75pt,black] (2.5,-2.5) to (2.5,2.5);
    \draw[line width=0.75pt,black] (2.5,2.5) to (-2.5,2.5);
    \draw[line width=0.75pt,black] (-2.5,2.5) to (-2.5,-2.5);
    \draw[line width=1.5pt,red] (-2.5,0) to (0,0);
    \draw[line width=0.75pt,black] (2.5,2.75) to (-2.5,2.75);
    \draw[->,line width=1.5pt,black] (0.0,3.0) to (0.0,3.5);
    \node[ ] at (0,3.85) {$\Delta\disp$};
    \draw[line width=0.75pt,black] (-2.5,-2.5) to (-2.75,-2.75);
    \draw[line width=0.75pt,black] (-2.25,-2.5) to (-2.5,-2.75);
    \draw[line width=0.75pt,black] (-2.0,-2.5) to (-2.25,-2.75);
    \draw[line width=0.75pt,black] (-1.75,-2.5) to (-2.0,-2.75);
    \draw[line width=0.75pt,black] (-1.5,-2.5) to (-1.75,-2.75);
    \draw[line width=0.75pt,black] (-1.25,-2.5) to (-1.5,-2.75);
    \draw[line width=0.75pt,black] (-1.0,-2.5) to (-1.25,-2.75);
    \draw[line width=0.75pt,black] (-0.75,-2.5) to (-1.0,-2.75);
    \draw[line width=0.75pt,black] (-0.5,-2.5) to (-0.75,-2.75);
    \draw[line width=0.75pt,black] (-0.25,-2.5) to (-0.5,-2.75);
    \draw[line width=0.75pt,black] (0.0,-2.5) to (-0.25,-2.75);
    \draw[line width=0.75pt,black] (0.25,-2.5) to (0.0,-2.75);
    \draw[line width=0.75pt,black] (0.5,-2.5) to (0.25,-2.75);
    \draw[line width=0.75pt,black] (0.75,-2.5) to (0.5,-2.75);
    \draw[line width=0.75pt,black] (1.0,-2.5) to (0.75,-2.75);
    \draw[line width=0.75pt,black] (1.25,-2.5) to (1.0,-2.75);
    \draw[line width=0.75pt,black] (1.5,-2.5) to (1.25,-2.75);
    \draw[line width=0.75pt,black] (1.75,-2.5) to (1.5,-2.75);
    \draw[line width=0.75pt,black] (2.0,-2.5) to (1.75,-2.75);
    \draw[line width=0.75pt,black] (2.25,-2.5) to (2.0,-2.75);
    \draw[line width=0.75pt,black] (2.5,-2.5) to (2.25,-2.75);
    \end{tikzpicture}
\caption{SENT experiment}
\label{sec4:fig:tension}
\end{minipage}
\begin{minipage}[b]{0.33\linewidth}
\centering
\begin{tikzpicture}[scale=0.7]
    \coordinate (K) at (0,0);
    \draw[line width=0.75pt,black] (-2.5,-2.5) to (2.5,-2.5);
    \draw[line width=0.75pt,black] (2.5,-2.5) to (2.5,2.5);
    \draw[line width=0.75pt,black] (2.5,2.5) to (-2.5,2.5);
    \draw[line width=0.75pt,black] (-2.5,2.5) to (-2.5,-2.5);
    \draw[line width=1.5pt,red] (-2.5,0) to (0,0);
    \draw[line width=0.75pt,black] (2.5,2.75) to (-2.5,2.75);
    \draw[->,line width=1.5pt,black] (0.1,3.25) to (0.8,3.25);
    \node[ ] at (-0.5,3.25) {$\Delta\disp$};
    \draw[line width=0.75pt,black] (-2.5,-2.5) to (-2.75,-2.75);
    \draw[line width=0.75pt,black] (-2.25,-2.5) to (-2.5,-2.75);
    \draw[line width=0.75pt,black] (-2.0,-2.5) to (-2.25,-2.75);
    \draw[line width=0.75pt,black] (-1.75,-2.5) to (-2.0,-2.75);
    \draw[line width=0.75pt,black] (-1.5,-2.5) to (-1.75,-2.75);
    \draw[line width=0.75pt,black] (-1.25,-2.5) to (-1.5,-2.75);
    \draw[line width=0.75pt,black] (-1.0,-2.5) to (-1.25,-2.75);
    \draw[line width=0.75pt,black] (-0.75,-2.5) to (-1.0,-2.75);
    \draw[line width=0.75pt,black] (-0.5,-2.5) to (-0.75,-2.75);
    \draw[line width=0.75pt,black] (-0.25,-2.5) to (-0.5,-2.75);
    \draw[line width=0.75pt,black] (0.0,-2.5) to (-0.25,-2.75);
    \draw[line width=0.75pt,black] (0.25,-2.5) to (0.0,-2.75);
    \draw[line width=0.75pt,black] (0.5,-2.5) to (0.25,-2.75);
    \draw[line width=0.75pt,black] (0.75,-2.5) to (0.5,-2.75);
    \draw[line width=0.75pt,black] (1.0,-2.5) to (0.75,-2.75);
    \draw[line width=0.75pt,black] (1.25,-2.5) to (1.0,-2.75);
    \draw[line width=0.75pt,black] (1.5,-2.5) to (1.25,-2.75);
    \draw[line width=0.75pt,black] (1.75,-2.5) to (1.5,-2.75);
    \draw[line width=0.75pt,black] (2.0,-2.5) to (1.75,-2.75);
    \draw[line width=0.75pt,black] (2.25,-2.5) to (2.0,-2.75);
    \draw[line width=0.75pt,black] (2.5,-2.5) to (2.25,-2.75);
    \draw[fill=gray!50] (-2.75,1.25) -- (-2.25,1.25) -- (-2.5,1.75)-- (-2.75,1.25);
    \draw[fill=black!75] (-2.5,1.05) circle (0.2);
    \draw[line width=1pt,black] (-2.25,0.8) to (-2.75,0.8);
    \draw[fill=gray!50] (-2.75,-1.5) -- (-2.25,-1.5) -- (-2.5,-1.)-- (-2.75,-1.5);
    \draw[fill=black!75] (-2.5,-1.7) circle (0.2);
    \draw[line width=1pt,black] (-2.25,-1.95) to (-2.75,-1.95);
    \draw[fill=gray!50] (2.75,1.25) -- (2.25,1.25) -- (2.5,1.75)-- (2.75,1.25);
    \draw[fill=black!75] (2.5,1.05) circle (0.2);
    \draw[line width=1pt,black] (2.25,0.8) to (2.75,0.8);
    \draw[fill=gray!50] (2.75,-1.5) -- (2.25,-1.5) -- (2.5,-1.)-- (2.75,-1.5);
    \draw[fill=black!75] (2.5,-1.7) circle (0.2);
    \draw[line width=1pt,black] (2.25,-1.95) to (2.75,-1.95);
    \end{tikzpicture}
\caption{SENS experiment}
\label{sec4:fig:miehe_shear}
\end{minipage}
\begin{minipage}[b]{0.33\linewidth}
\centering
\begin{tabular}{ll} \hline
  \textbf{Parameters} & \textbf{Value} \\ \hline
  Model & Brittle-AT2 \\
  $\lambda$ & 121.154 [GPa] \\
  $\mu$ & 80.769 [GPa] \\
  $\gc$ & 2700 [N/m] \\
  $\l$ & 1.5e-2 [mm] \\
  Element size & $\l/2$ \\ \hline
  \end{tabular}
\captionof{table}{Parameters}
\label{sec4:table:miehe_tension_shear}
\end{minipage}
\end{figure}


\begin{figure}[!ht]
  \begin{subfigure}[t]{0.45\textwidth}
  \centering
    \begin{tikzpicture}[thick,scale=0.95, every node/.style={scale=1.175}]
    \begin{axis}[legend style={draw=none}, legend columns = 3,
      transpose legend, ylabel={Load\:[kN]},xlabel={Displacement\:[mm]}, xmin=0, ymin=0, xmax=0.0065, ymax=1.1, yticklabel style={/pgf/number format/.cd,fixed,precision=2},
                 every axis plot/.append style={very thick}]
    \pgfplotstableread{./Data/SingleEdgeNotchedTension/NotchTension_Miehe.txt}\Ddata;
    \pgfplotstableread{./Data/SingleEdgeNotchedTension/NotchTension_Ambati.txt}\Cdata;
    \pgfplotstableread{./Data/SingleEdgeNotchedTension/NotchTension_Pen_Lodi.txt}\Bdata;
    \pgfplotstableread{./Data/SingleEdgeNotchedTension/NotchTension_Aug_Lodi.txt}\Adata;
    \addplot [ 
           color=black, 
           mark=*, 
           mark size=0.75pt, 
         ]
         table
         [
           x expr=\thisrowno{2}, 
           y expr=\thisrowno{1}
         ] {\Adata};
         \addlegendentry{LMM}
    \addplot [ 
           color=red, 
           mark=*, 
           mark size=0.25pt, 
         ]
         table
         [
           x expr=\thisrowno{2}, 
           y expr=\thisrowno{1}
         ] {\Bdata};
         \addlegendentry{Pen.}
    \addplot [ 
           color=blue, 
           mark=*, 
           mark size=0.25pt, 
         ]
         table
         [
           x expr=\thisrowno{0}, 
           y expr=\thisrowno{1}
         ] {\Cdata};
         \addlegendentry{Ambati \cite{Ambati2014}}
    \addplot [ 
           color=green, 
           mark=*, 
           mark size=0.25pt, 
         ]
         table
         [
           x expr=\thisrowno{0}, 
           y expr=\thisrowno{1}
         ] {\Ddata};
         \addlegendentry{Miehe \cite{Miehe2010a}}
    \end{axis}
    \end{tikzpicture}
    \caption{ }
    \label{sec4:fig:miehe_tension_lodi}
  \end{subfigure}
  \hfill
  \begin{subfigure}[t]{0.45\textwidth}
  \centering
    \begin{tikzpicture}
    \node[inner sep=0pt] () at (0,0)
    {\includegraphics[width=5.5cm,trim=5cm 0cm 5cm 0cm, clip]{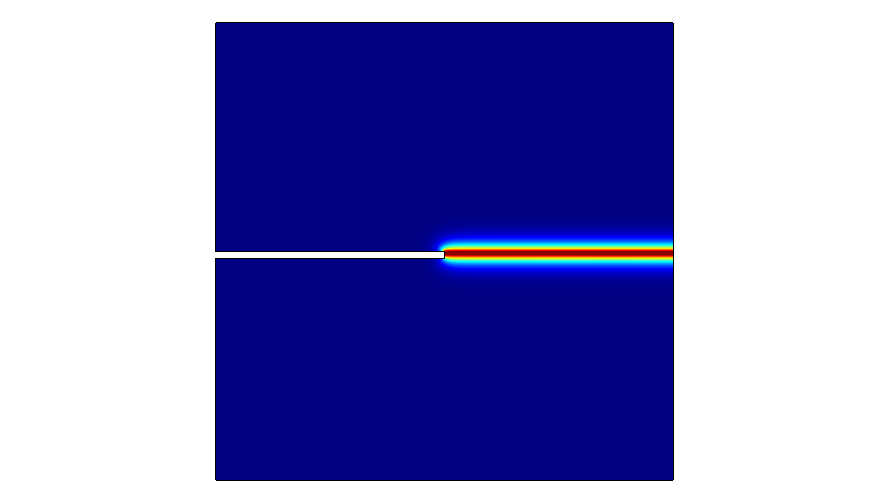}};
    \node[inner sep=0pt] () at (-1.15,-3.05)
    {\begin{axis}[
    hide axis,
    scale only axis,
    height=0pt,
    width=0pt,
    colormap/bluered,
    colorbar horizontal,
    point meta min=0,
    point meta max=1,
    colorbar style={
        width=4.75cm,
        xtick={0,0.5,1.0},
        xticklabel style = {yshift=-0.075cm}
    }]
    \addplot [draw=none] coordinates {(0,0)};
    \end{axis}};
    \node[inner sep=0pt] () at (0,-3.75) {$\pf$};
    \end{tikzpicture}
    \caption{ }
    \label{sec4:fig:pf_sen_tension_failure}
  \end{subfigure}
  \caption{Figure (a) presents the load-displacement curves for the single edge notched specimen under tension. Here, LMM and Pen. refer to the Lagrange Multiplier Method and the Penalty method respectively. Figure (b) shows the distribution of the phase-field variable at the final step of the analysis.}
\end{figure}

Figure \ref{sec4:fig:miehe_tension_lodi} presents the load-displacement curves obtained using the Lagrange Multiplier Method (LMM) and Penalty (Pen.) method, along with those from the literature \cite{Miehe2010a,Ambati2014}. Both methods yield a similar behaviour compared to \cite{Miehe2010a}, in terms of the peak load and the post-peak behaviour. Moreover, the phase-field fracture topology at failure in Figure \ref{sec4:fig:pf_sen_tension_failure} is also identical to those reported in the aforementioned literature.

\subsection{Brittle: Single Edge Notched specimen under Shear (SENS)}\label{sec5:SENshear}

In order to perform a shear test, the SEN specimen is loaded horizontally along the top edge as shown in Figure \ref{sec4:fig:miehe_shear}. The material properties remain same as presented in Table \ref{sec4:table:miehe_tension_shear}. A quasi-static loading is applied to the top boundary in the form of prescribed displacement increment $\Delta\disp = 1e-4$[mm] for the first $90$ steps, following which it is changed to $1e-5$[mm]. Furthermore, the bottom boundary remains fixed, and roller support is implemented in left and right edges thereby restricting the vertical displacement.

\begin{figure}[!ht]
\centering
  \begin{subfigure}[t]{0.45\textwidth}
  \centering
    \begin{tikzpicture}[thick,scale=0.95, every node/.style={scale=1.175}]
    \begin{axis}[legend style={draw=none}, legend columns = 2,
      transpose legend, ylabel={Load\:[kN]},xlabel={Displacement\:[mm]}, xmin=0, ymin=0, xmax=0.015, ymax=1., yticklabel style={/pgf/number format/.cd,fixed,precision=2},
                 every axis plot/.append style={very thick}]
    \pgfplotstableread{./Data/SingleEdgeNotchedShear/NotchShear_Miehe.txt}\Ddata;
    \pgfplotstableread{./Data/SingleEdgeNotchedShear/NotchShear_Ambati.txt}\Cdata;
    \pgfplotstableread{./Data/SingleEdgeNotchedShear/NotchShear_Pen.txt}\Bdata;
    \pgfplotstableread{./Data/SingleEdgeNotchedShear/NotchShear_Lodi.txt}\Adata;
    \addplot [ 
           color=black, 
           mark=*, 
           mark size=0.75pt, 
         ]
         table
         [
           x expr=\thisrowno{2}, 
           y expr=\thisrowno{1}
         ] {\Adata};
         \addlegendentry{LMM}
    \addplot [ 
           color=red, 
           mark=*, 
           mark size=0.25pt, 
         ]
         table
         [
           x expr=\thisrowno{2}, 
           y expr=\thisrowno{1}
         ] {\Bdata};
         \addlegendentry{Pen.}     
    \addplot [ 
           color=green, 
           mark=*, 
           mark size=0.25pt, 
         ]
         table
         [
           x expr=\thisrowno{0}, 
           y expr=\thisrowno{1}
         ] {\Ddata};
         \addlegendentry{Miehe \cite{Miehe2010a}}
    \end{axis}
    \end{tikzpicture}
    \caption{ }
    \label{sec4:fig:miehe_shear_lodi}
  \end{subfigure}
  \hfill
  \begin{subfigure}[t]{0.45\textwidth}
  \centering
    \begin{tikzpicture}
    \node[inner sep=0pt] () at (0,0)
    {\includegraphics[width=5.5cm,trim=5cm 0cm 5cm 0cm, clip]{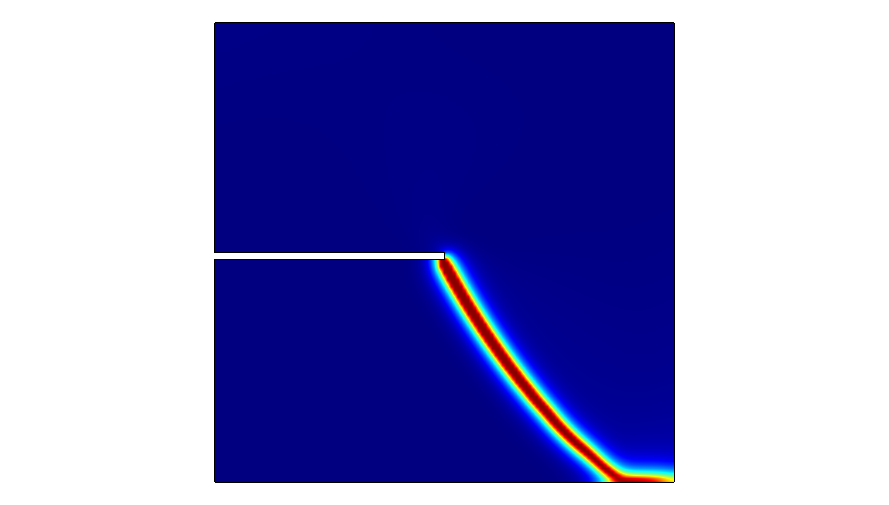}};
    \node[inner sep=0pt] () at (-1.15,-3.05)
    {\begin{axis}[
    hide axis,
    scale only axis,
    height=0pt,
    width=0pt,
    colormap/bluered,
    colorbar horizontal,
    point meta min=0,
    point meta max=1,
    colorbar style={
        width=4.75cm,
        xtick={0,0.5,1.0},
        xticklabel style = {yshift=-0.075cm}
    }]
    \addplot [draw=none] coordinates {(0,0)};
    \end{axis}};
    \node[inner sep=0pt] () at (0,-3.75) {$\pf$};
    \end{tikzpicture}
    \caption{ }
    \label{sec4:fig:pf_sen_shear_failure}
  \end{subfigure}
  \caption{Figure (a) presents the load-displacement curves for the single edge notched specimen under shear. Here, LMM and Pen. refer to the Lagrange Multiplier Method and the Penalty method respectively. Figure (b) shows the distribution of the phase-field variable at the final step of the analysis.}
  \label{sec4:fig:miehe_shear_lodi_iter}
\end{figure}

Figure \ref{sec4:fig:miehe_shear_lodi} presents the load-displacement curves obtained using the Lagrange Multiplier Method and the Penalty method, along with those from the literature \cite{Miehe2010a}. The former methods result in ~$9\%$ higher peak load estimation as compared to \cite{Miehe2010a}. However, the phase-field fracture topology in the final step of the analysis is consistent with those reported in the aforementioned literature.

\subsection{Brittle: Notched concrete specimen with a hole}\label{sec5:NotchHole}

\begin{figure}[ht]
\begin{minipage}[b]{0.33\linewidth}
\centering
\begin{tikzpicture}[scale=0.5]
    \coordinate (K) at (0,0);
    \draw[line width=0.75pt,black] (0,0) to (6.5,0);
    \draw[line width=0.75pt,black] (6.5,0) to (6.5,12.0);
    \draw[line width=0.75pt,black] (6.5,12.0) to (0.0,12.0);
    \draw[line width=0.75pt,black] (0.0,12.0) to (0.0,0.0);
    \draw[line width=1.5pt,red] (0.0,6.5) to (1,6.5);
    \draw[fill=black!25] (2,2) circle (0.5);
    \draw[line width=0.75pt,black] (3.65,5.1) circle (1);
    \draw[fill=black!25] (2,10) circle (0.5);
    \draw[->,line width=1.5pt,black] (2,10.65) to (2,11.05);
    \node[] at (2,11.5) {$\Delta\disp$};
    \end{tikzpicture}
\caption{Geometry and constraints}
\label{sec4:fig:notched_hole}
\end{minipage}
\begin{minipage}[b]{0.33\linewidth}
\centering
    \includegraphics[trim=0cm 0cm 0cm 0cm,clip,width=3.25cm]{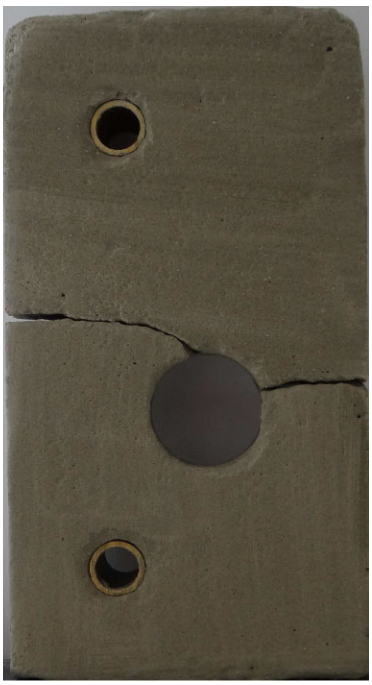}
\caption{Experimental results \cite{Ambati2014}}
\label{sec4:fig:notched_hole_exp}
\end{minipage}
\begin{minipage}[b]{0.33\linewidth}
\centering
\begin{tabular}{ll} \hline
  \textbf{Symbol}  & \textbf{Value} \\ \hline
  Model & Brittle-AT2 \\
  $\lambda$  & 1.94 [kN/mm$^2$] \\
  $\mu$ & 2.45 [kN/mm$^2$] \\
  $\gc$ & 2.28e-3 [kN/mm] \\
  $\l$ & 0.25 [mm] \\
  Element size & $\l/2$ \\ \hline
  \end{tabular}
\captionof{table}{Material properties \cite{ma14081913} }
\label{sec4:table:notched_hole}
\end{minipage}
\end{figure}

A notched concrete specimen with a hole, shown in Figure \ref{sec4:fig:notched_hole} is considered in this study. The experimental and numerical analysis of the same has been carried out earlier in \cite{Ambati2014}. The specimen has dimensions 65 x 120 [mm$^2$], the hole being 20 [mm] in diameter located at (36.5 [mm], 51[mm]). Moreover, a notch, 10 [mm] in length is located at 65 [mm] height from the bottom of the plate. The model parameters used in the simulation is presented in Table \ref{sec4:table:notched_hole}, while the experimentally observed fracture pattern is shown in Figure \ref{sec4:fig:notched_hole_exp}. As shown in Figure \ref{sec4:fig:notched_hole}, the plate is loaded via the upper pin (in grey) with displacement increment $\Delta\disp = 1e-3$ mm, while the lower pin (in grey) is remains fixed.

\begin{figure}[!ht]
\centering
  \begin{subfigure}[t]{0.4\textwidth}
  \centering
    \begin{tikzpicture}[thick,scale=0.95, every node/.style={scale=1.155}]
    \begin{axis}[legend style={draw=none}, legend columns = 3,
      transpose legend, ylabel={Load\:[kN]},xlabel={Displacement\:[mm]}, xmin=0, ymin=0, xmax=2, ymax=1.0, yticklabel style={/pgf/number format/.cd,fixed,precision=2},
                 every axis plot/.append style={very thick}]
    \pgfplotstableread{./Data/NotchedHole/NotchHole_Ambati_Lodi.txt}\Cdata;
    \pgfplotstableread{./Data/NotchedHole/NotchHole_Pen_Lodi_l0.25.txt}\Bdata;
    \pgfplotstableread{./Data/NotchedHole/NotchHole_Aug_Lodi_l0.25.txt}\Adata;
    \pgfplotstableread{./Data/NotchedHole/NotchedHole_Paneda.txt}\Ddata;
    \addplot [ 
           color=black, 
           mark=*, 
           mark size=0.75pt, 
         ]
         table
         [
           x expr=\thisrowno{2}, 
           y expr=\thisrowno{1}
         ] {\Adata};
         \addlegendentry{LMM}
    \addplot [ 
           color=red, 
           mark=*, 
           mark size=0.25pt, 
         ]
         table
         [
           x expr=\thisrowno{2}, 
           y expr=\thisrowno{1}
         ] {\Bdata};
         \addlegendentry{Pen.}     
    \addplot [ 
           color=blue, 
           mark=*, 
           mark size=0.25pt, 
         ]
         table
         [
           x expr=\thisrowno{0}, 
           y expr=\thisrowno{1}
         ] {\Cdata};
         \addlegendentry{Ambati \cite{Ambati2014}}
    \addplot [ 
           color=green, 
           mark=*, 
           mark size=0.25pt, 
         ]
         table
         [
           x expr=\thisrowno{0}, 
           y expr=\thisrowno{1}
         ] {\Ddata};
         \addlegendentry{Navidtehrani \cite{ma14081913}}
    \end{axis}
    \end{tikzpicture}
    \caption{ }
    \label{sec4:fig:notch_hole_lodi}
  \end{subfigure}
  \hfill
  \begin{subfigure}[t]{0.45\textwidth}
  \centering
    \begin{tikzpicture}
    \node[inner sep=0pt] () at (0,0)
    {\includegraphics[width=7cm,trim=7cm 3.5cm 7cm 3.5cm, clip]{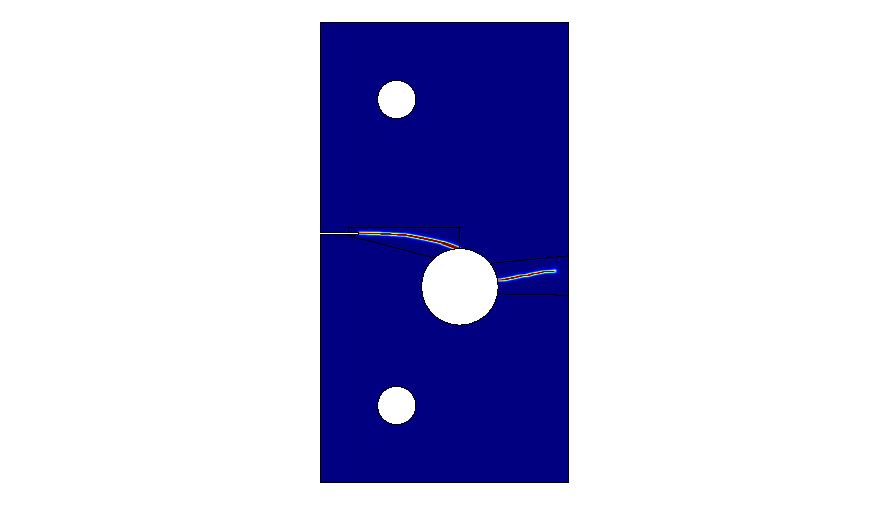}};
    \node[inner sep=0pt] () at (-1.2,-3.0)
    {\begin{axis}[
    hide axis,
    scale only axis,
    height=0pt,
    width=0pt,
    colormap/bluered,
    colorbar horizontal,
    point meta min=0,
    point meta max=1,
    colorbar style={
        width=4.75cm,
        xtick={0,0.5,1.0},
        xticklabel style = {yshift=-0.075cm}
    }]
    \addplot [draw=none] coordinates {(0,0)};
    \end{axis}};
    \node[inner sep=0pt] () at (0,-3.6) {$\pf$};
    \end{tikzpicture}
    \caption{ }
    \label{sec4:fig:pf_notch_hole_failure}
  \end{subfigure}
  \caption{Figure (a) presents the load-displacement curves for the notched specimen with hole test, obtained from the current implementation, \cite{Ambati2014}, and \cite{ma14081913}. Here, LMM and Pen. refer to the Lagrange Multiplier Method and the Penalty method respectively. Figure (b) shows the distribution of the phase-field variable at the final step of the analysis in a section of the specimen.}
  \label{sec4:fig:notch_hole_lodi_pf}
\end{figure}

Figure \ref{sec4:fig:notch_hole_lodi} presents the load-displacement curves obtained using the Lagrange Multiplier Method and the Penalty method, along with those from the literature \cite{Ambati2014,ma14081913}. The curves corresponding to the Lagrange Multiplier Method and the Penalty method coincide, and predict peak loads close to those obtained in \cite{Ambati2014} and \cite{ma14081913}. However, the load-displacement curves are significantly different from \cite{Ambati2014}. The reason for this behaviour is not investigated in this manuscript. However, it is important to note that in this work, the upper and lower pins are assumed to be rigid connectors. As such, the loading or fixed boundary conditions is applied w.r.t. the centre of the pin. Furthermore, \cite{Ambati2014} opted for $\l = 0.1$ [mm], in this manuscript as well as in \cite{ma14081913}, $\l = 0.25$ [mm]. Finally, from Figure \ref{sec4:fig:pf_notch_hole_failure}, it observed that the final phase-field fracture topology does match the experimentally observed fracture patterns in Figure \ref{sec4:fig:notched_hole_exp}.

\subsection{Quasi-brittle: Concrete three-point bending}\label{sec5:concreteTPB}

A three-point bending experiment on a notched concrete beam reported in \cite{Rots1988} is considered here. The beam has dimensions $450\times100$ [mm$^2$], and has a notch $5\times50$ [mm$^2$]. A schematic of the beam along with the loading conditions is presented in Figure \ref{sec4:fig:3pt_bending}. Displacement-based load increments of $\Delta\disp = 1e-3$ [mm] is enforced throughout the simulation. The model parameters are presented in Table \ref{sec4:table:3pt_bending}.

\begin{figure}[ht]
\begin{minipage}[b]{0.45\linewidth}
\centering
\begin{tikzpicture}[scale=0.7]
    \coordinate (K) at (0,0);
    \draw[line width=0.75pt,black] (4,2) to (-4,2);
    \draw[line width=0.75pt,black] (4,2) to (4,0);
    \draw[line width=0.75pt,black] (4,0) to (0.1,0);
    \draw[line width=0.75pt,black] (0.05,0) to (0.05,1.0);
    \draw[line width=0.75pt,black] (0.05,1.0) to (-0.05,1.0);
    \draw[line width=0.75pt,black] (-0.05,1.0) to (-0.05,0.0);
    \draw[line width=0.75pt,black] (-0.05,0.0) to (-4.0,0.0);
    \draw[line width=0.75pt,black] (-4.0,0.0) to (-4.0,2);
    \draw[->,line width=1.5pt,black] (0.0,2.5) to (0.0,2.);
    \node[ ] at (0,2.75) {$\Delta\disp$};
    \draw[fill=gray!50] (-4.25,-0.5) -- (-3.75,-0.5) -- (-4,0)-- (-4.25,-0.5);
    \draw[line width=1pt,black] (-4.5,-0.5) to (-3.5,-0.5);
    \draw[line width=1pt,black] (-4.5,-0.5) to (-4.75,-0.75);
    \draw[line width=1pt,black] (-4.25,-0.5) to (-4.5,-0.75);
    \draw[line width=1pt,black] (-4,-0.5) to (-4.25,-0.75);
    \draw[line width=1pt,black] (-3.75,-0.5) to (-4,-0.75);
    \draw[line width=1pt,black] (-3.5,-0.5) to (-3.75,-0.75);
    \draw[fill=gray!50] (4.25,-0.5) -- (3.75,-0.5) -- (4,0)-- (4.25,-0.5);
    \draw[fill=black!75] (4,-0.7) circle (0.2);
    \draw[line width=1pt,black] (4.5,-0.95) to (3.5,-0.95);
    \end{tikzpicture}
\caption{Three point bending test}
\label{sec4:fig:3pt_bending}
\end{minipage}
\hspace{0.5cm}
\begin{minipage}[b]{0.45\linewidth}
\centering
\begin{tabular}{ll} \hline
  \textbf{Parameters}  & \textbf{Value} \\ \hline
  Model & Quasi-brittle \\
  Softening & Cornellisen et. al. \cite{cornelissen1986experimental} \\
  $E_0$  & 2e4 [MPa] \\
  $\nu$  & 0.2 [-] \\
  $f_t$  & 2.4 [MPa] \\
  $\gc$  & 0.113 [N/mm] \\
  $\l$  & 2.5 [mm] \\
  Element size & $\l/5$ \\ \hline
  \end{tabular}
\captionof{table}{Parameters for three point bending test}
\label{sec4:table:3pt_bending}
\end{minipage}
\end{figure}

\begin{figure}[!ht]
\centering
  \begin{subfigure}[t]{0.45\textwidth}
  \centering
    \begin{tikzpicture}[thick,scale=0.95, every node/.style={scale=1.175}]
    \begin{axis}[legend style={draw=none}, legend columns = 3,
      transpose legend, ylabel={Load\:[kN]},xlabel={Displacement\:[mm]}, xmin=0, ymin=0, xmax=1.0, ymax=1.75, yticklabel style={
        /pgf/number format/fixed,
        /pgf/number format/precision=5
        },
        scaled y ticks=false,
        xticklabel style={
        /pgf/number format/fixed,
        /pgf/number format/precision=5
        },
        scaled x ticks=false,
        every axis plot/.append style={very thick}]
    \pgfplotstableread{./Data/RotsTPB/MM_RotsTPB_Lodi.txt}\Ddata;
    \pgfplotstableread{./Data/RotsTPB/Pen_RotsTPB_Lodi.txt}\Cdata;
    \pgfplotstableread[col sep = comma]{./Data/RotsTPB/RotsTPBUpperLodi.txt}\Bdata;
    \pgfplotstableread[col sep = comma]{./Data/RotsTPB/RotsTPBLowerLodi.txt}\Adata;
    \addplot [name path = A, 
           color=green!20, 
           mark=*, 
           mark size=0.25pt, 
           forget plot
         ]
         table
         [
           x expr=\thisrowno{0}, 
           y expr=\thisrowno{1}
         ] {\Adata};
    \addplot [name path = B,
           color=green!20, 
           mark=*, 
           mark size=0.25pt, 
           forget plot
         ]
         table
         [
           x expr=\thisrowno{0}, 
           y expr=\thisrowno{1}
         ] {\Bdata};
    \addplot [green!20, forget plot] fill between [of = A and B]; 
    \addplot [ 
           color=black, 
           mark=*, 
           mark size=0.75pt, 
         ]
         table
         [
           x expr=\thisrowno{2}, 
           y expr=\thisrowno{1}
         ] {\Ddata};
         \addlegendentry{LMM}
    \addplot [ 
           color=red, 
           mark=*, 
           mark size=0.25pt, 
         ]
         table
         [
           x expr=\thisrowno{2}, 
           y expr=\thisrowno{1}
         ] {\Cdata};
         \addlegendentry{Pen.}
    \end{axis}
    \end{tikzpicture}
    \caption{Load-displacement plot}
    \label{sec4:fig:3pt_bending_lodi}
  \end{subfigure}
  \hspace{5mm}
  \begin{subfigure}[t]{0.45\textwidth}
  \centering
    \begin{tikzpicture}
    \node[inner sep=0pt] () at (0,0)
    {\includegraphics[width=6cm,trim=5cm 0cm 5cm 0cm, clip]{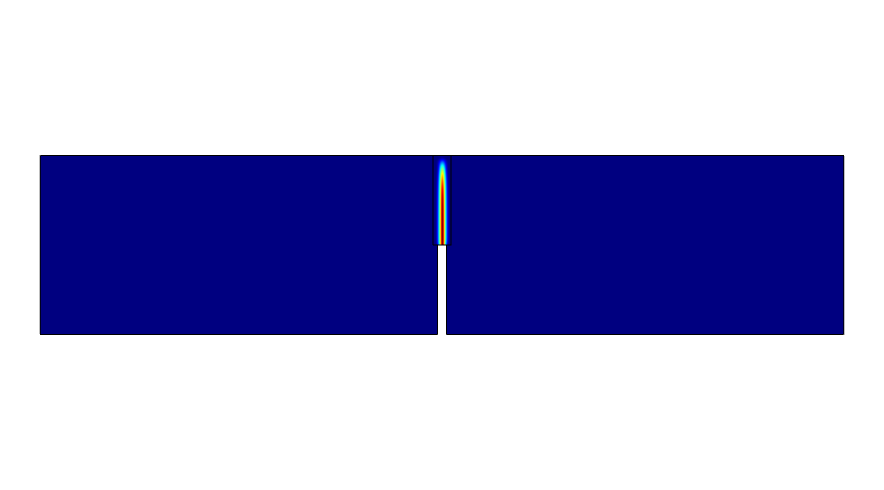}};
    \node[inner sep=0pt] () at (-1.2,-1.75)
    {\begin{axis}[
    hide axis,
    scale only axis,
    height=0pt,
    width=0pt,
    colormap/bluered,
    colorbar horizontal,
    point meta min=0,
    point meta max=1,
    colorbar style={
        width=4.75cm,
        xtick={0,0.5,1.0},
        xticklabel style = {yshift=-0.075cm}
    }]
    \addplot [draw=none] coordinates {(0,0)};
    \end{axis}};
    \node[inner sep=0pt] () at (0,-2.45) {$\pf$};
    \end{tikzpicture}
    \caption{ }
    \label{sec4:fig:pf_3pt_bending_failure}
  \end{subfigure}
  \caption{Figure (a) presents the load-displacement curves for the concrete three-point bending test. Here, LMM and Pen. refer to the Lagrange Multiplier Method and the Penalty method respectively. The experimental range is represented by the shaded area. Figure (b) shows the distribution of the phase-field variable at the final step of the analysis in a section of the beam.}
  \label{sec4:fig:3pt_bending_results}
\end{figure}

Figure \ref{sec4:fig:3pt_bending_lodi} presents the load-displacement curves obtained using  the Lagrange Multiplier Method and the Penalty method, along with those obtained from the experiments \cite{Rots1988}. It is observed that both methods yield identical curves. These curves are close to the experimental upper bound (shaded region). Furthermore, the phase-field fracture topology at the final step in the analysis in Figure \ref{sec4:fig:pf_3pt_bending_failure} is also identical to that reported in \cite{Rots1988}. 

\subsection{Quasi-brittle: Winkler L-panel}\label{sec5:Lpanel}

The L-shaped panel studied by \cite{winkler2001traglastuntersuchungen,unger2007modelling} is considered in this sub-section. Figure \ref{sec4:fig:Lpanel} shows the geometry as well as the loading conditions. The longer edges of the panel are 500 [mm] and the smaller edges are 250 [mm]. The loading is applied on the edge marked in blue, 30 [mm] in length, and is in the form of displacement-based increments of $\Delta\disp = 1e-3$ [mm]. The additional model parameters are presented in Table \ref{sec4:table:Lpanel_params}.

\begin{figure}[ht]
\begin{minipage}[b]{0.45\linewidth}
\centering
\begin{tikzpicture}[scale=0.6]
    \coordinate (K) at (0,0);
    \draw[line width=0.75pt,black] (-4,-4) to (0,-4);
    \draw[line width=0.75pt,black] (0,-4) to (0,0);
    \draw[line width=0.75pt,black] (0,0) to (4,0);
    \draw[line width=0.75pt,black] (4,0)  to (4,4);
    \draw[line width=0.75pt,black] (4,4) to (-4,4);
    \draw[line width=0.75pt,black] (-4,4) to (-4,-4);
    \draw[line width=2.0pt,blue] (3.5,0) to (4,0);
    \draw[->,line width=1.5pt,black] (3.75,-0.9) to (3.75,-0.25);
    \node[ ] at (3.75,-1.15) {$\Delta\disp$};
    \draw[line width=1pt,black] (0,-4) to (-0.25,-4.25);
    \draw[line width=1pt,black] (-0.25,-4.) to (-0.5,-4.25);
    \draw[line width=1pt,black] (-0.5,-4.) to (-0.75,-4.25);
    \draw[line width=1pt,black] (-0.75,-4.) to (-1,-4.25);
    \draw[line width=1pt,black] (-1,-4.) to (-1.25,-4.25);
    \draw[line width=1pt,black] (-1.25,-4.) to (-1.5,-4.25);
    \draw[line width=1pt,black] (-1.5,-4.) to (-1.75,-4.25);
    \draw[line width=1pt,black] (-1.75,-4.) to (-2.,-4.25);
    \draw[line width=1pt,black] (-2.,-4.) to (-2.25,-4.25);
    \draw[line width=1pt,black] (-2.25,-4.) to (-2.5,-4.25);
    \draw[line width=1pt,black] (-2.5,-4.) to (-2.75,-4.25);
    \draw[line width=1pt,black] (-2.75,-4.) to (-3.,-4.25);
    \draw[line width=1pt,black] (-3.,-4.) to (-3.25,-4.25);
    \draw[line width=1pt,black] (-3.25,-4.) to (-3.5,-4.25);
    \draw[line width=1pt,black] (-3.5,-4.) to (-3.75,-4.25);
    \draw[line width=1pt,black] (-3.75,-4.) to (-4.,-4.25);
    \draw[line width=1pt,black] (-4.,-4.) to (-4.25,-4.25);
    \end{tikzpicture}
\caption{Winkler L-panel}
\label{sec4:fig:Lpanel}
\end{minipage}
\hspace{0.5cm}
\begin{minipage}[b]{0.45\linewidth}
\centering
\begin{tabular}{ll} \hline
  \textbf{Parameters}  & \textbf{Value} \\ \hline
  Model & Quasi-brittle \\
  Softening & Cornellisen et. al. \cite{cornelissen1986experimental} \\
  $E_0$  & 2.0e4 [MPa] \\
  $\nu$  & 0.18 [-] \\
  $f_t$  & 2.5 [MPa] \\
  $\gc$  & 0.130 [N/mm] \\
  $\l$  & 5 [mm] \\
  Element size & $\l/5$ \\ \hline
  \end{tabular}
\captionof{table}{Parameters for L-shaped panel test \cite{unger2007modelling}}
\label{sec4:table:Lpanel_params}
\end{minipage}
\end{figure}

\begin{figure}[!ht]
\centering
  \begin{subfigure}[t]{0.45\textwidth}
  \centering
    \begin{tikzpicture}[thick,scale=0.95, every node/.style={scale=1.175}]
    \begin{axis}[legend style={draw=none}, legend columns = 3,
      transpose legend, ylabel={Load\:[kN]},xlabel={Displacement\:[mm]}, xmin=0, ymin=0, xmax=1.0, ymax=8.0, yticklabel style={
        /pgf/number format/fixed,
        /pgf/number format/precision=5
        },
        scaled y ticks=false,
        xticklabel style={
        /pgf/number format/fixed,
        /pgf/number format/precision=5
        },
        scaled x ticks=false,
        every axis plot/.append style={very thick}]
    \pgfplotstableread{./Data/WinklerL/Pen_WinklerL_Lodi.txt}\Ddata;
    \pgfplotstableread{./Data/WinklerL/Pen_WinklerL_Lodi.txt}\Cdata;
    \pgfplotstableread[col sep = comma]{./Data/WinklerL/LPanelUpperLodi.txt}\Bdata;
    \pgfplotstableread[col sep = comma]{./Data/WinklerL/LPanelLowerLodi.txt}\Adata;
    \addplot [ name path = A,
           color=green!20, 
           mark=*, 
           mark size=0.01pt,
           forget plot
         ]
         table
         [
           x expr=\thisrowno{0}, 
           y expr=\thisrowno{1}
         ] {\Adata};
    \addplot [ name path = B,
           color=green!20, 
           mark=*, 
           mark size=0.01pt, 
           forget plot
         ]
         table
         [
           x expr=\thisrowno{0}, 
           y expr=\thisrowno{1}
         ] {\Bdata};
    \addplot [green!20, forget plot] fill between [of = A and B];     
    \addplot [ 
           color=black, 
           mark=*, 
           mark size=0.75pt, 
         ]
         table
         [
           x expr=\thisrowno{2}, 
           y expr=\thisrowno{1}
         ] {\Ddata};
         \addlegendentry{LMM}
    \addplot [ 
           color=red, 
           mark=*, 
           mark size=0.25pt, 
         ]
         table
         [
           x expr=\thisrowno{2}, 
           y expr=\thisrowno{1}
         ] {\Cdata};
         \addlegendentry{Pen.}
    \end{axis}
    \end{tikzpicture}
    \caption{Load-displacement plot}
    \label{sec4:fig:LPanel_lodi}
  \end{subfigure}
  \hspace{5mm}
  \begin{subfigure}[t]{0.45\textwidth}
  \centering
    \begin{tikzpicture}
    \node[inner sep=0pt] () at (0,1.5)
    {\includegraphics[width=10cm,trim=0cm 0cm 0cm 0cm, clip]{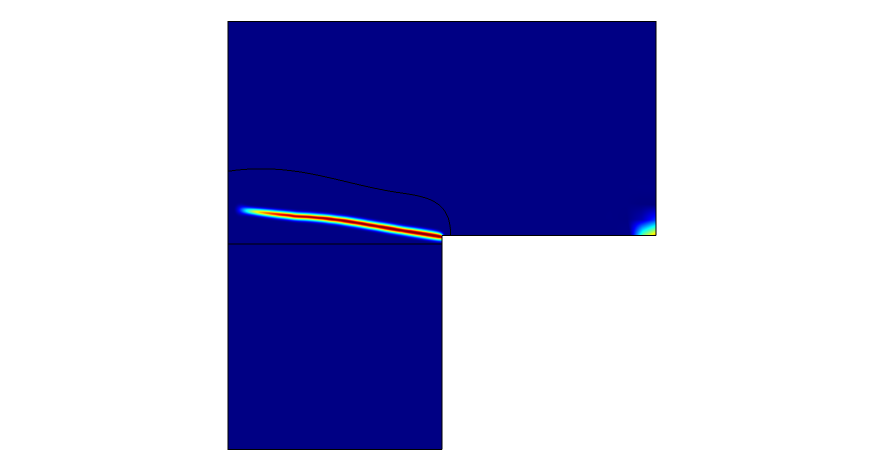}};
    \node[inner sep=0pt] () at (-1.2,-1.75)
    {\begin{axis}[
    hide axis,
    scale only axis,
    height=0pt,
    width=0pt,
    colormap/bluered,
    colorbar horizontal,
    point meta min=0,
    point meta max=1,
    colorbar style={
        width=4.75cm,
        xtick={0,0.5,1.0},
        xticklabel style = {yshift=-0.075cm}
    }]
    \addplot [draw=none] coordinates {(0,0)};
    \end{axis}};
    \node[inner sep=0pt] () at (0,-2.45) {$\pf$};
    \end{tikzpicture}
    \caption{ }
    \label{sec4:fig:pf_LPanel_failure}
  \end{subfigure}
  \caption{Figure (a) presents the load-displacement curves for the concrete Winkler L-panel test. Here, LMM and Pen. refer to the Lagrange Multiplier Method and the Penalty method respectively. The experimental range is represented by the shaded area. Figure (b) shows the distribution of the phase-field variable at the final step of the analysis in a section of the beam.}
  \label{sec4:fig:LPanel_results}
\end{figure}

Figure \ref{sec4:fig:LPanel_lodi} presents the load-displacement curves obtained using the Lagrange Multiplier Method and the Penalty method, along with those obtained from the experiments \cite{winkler2001traglastuntersuchungen,unger2007modelling}. It is observed that the curves obtained using both methods are identical. These curves are also close to the experimental range (shaded region). Furthermore, the phase-field fracture topology at the final step in the analysis in Figure \ref{sec4:fig:pf_LPanel_failure} is also identical to that reported in \cite{winkler2001traglastuntersuchungen,unger2007modelling}.

\section{Concluding Remarks}\label{sec6}

An alternative treatment of the phase-field fracture irreversibility constraint is presented in this manuscript. The fracture irreversibility constraint $h(\pf) \geq 0$ is transformed into an equivalent equality constraint. This is carried out upon introducing a slack variable $\theta$, and equating $h(\pf)$ to $\theta^2$. By definition, $\theta^2$ admits zero or positive values, thereby ensuring the fracture irreversibility constraint. The constraint is then introduced in the phase-field fracture energy functional using the Lagrange Multiplier Method and the Penalty method. The slack variable approach preserves the variational nature of the phase-field fracture problem, unlike the history variable approach in \cite{Miehe2010}. Furthermore, numerical experiments are conducted on benchmark brittle (AT2) and quasi-brittle fracture problems to demonstrate the efficacy of the proposed methods.

The scope for future studies include the extension towards complex multi-physics problems, large-scale problems, that maybe require special preconditioning techniques. Other future research directions may be towards accelerating the convergence rate using Anderson acceleration \cite{anderson1965iterative} or Large time increment method \cite{givoli2017latin,Bharali2017}.

\section{Software Implementation and Data Availability}

The numerical study in Section \ref{sec5} is carried out using the equation-based modelling approach in the software package COMSOL Multiphysics 5.6. The source files would be made available in Github repository of the corresponding author (\url{https://github.com/rbharali/PhaseFieldFractureCOMSOL/}).

\section*{Acknowledgements}

The financial support from the Swedish Research Council for Sustainable Development (FORMAS) under Grant 2018-01249 and the Swedish Research Council (VR) under Grant 2017-05192 is gratefully acknowledged. The authors would also like to thank Prof. Laura de Lorenzis (ETH Z\"urich) for insightful comments and suggestions.

\printbibliography

\end{document}